\documentclass[12pt]{article}

\usepackage{amsmath,amsthm,amstext,amscd,amssymb,euscript}
\usepackage{calrsfs}
\usepackage{epsfig}

\newcommand{\N}{{\mathbb N}}
\newcommand{\Z}{{\mathbb Z}}
\newcommand{\F}{{\mathcal F}}
\newcommand{\R}{{\mathbb R}}
\newcommand{\E}{{\mathbb E}}
\newcommand{\pee}{{\mathbb P}}
\newcommand{\co}{{\hat{\pee}}}
\newcommand{\eps}{\epsilon}
\newcommand{\var}{\textup{Var}}

\font\gfont=cmmi10 scaled \magstep{2.5}     
\newcommand{\gtau}{\hbox{\gfont \char28}}
\newcommand{\T}{{\gtau}}

\newcommand{\xinf}{x_{-\infty}^i}

\newcommand{\Xinf}{X_{-\infty}^i}
\newcommand{\Xinfi}{X_{-\infty}^{i-1}}

\newcommand{\yinf}{y_{-\infty}^i}

\newcommand{\Finf}{\F_{-\infty}^i}
\newcommand{\Finfi}{\F_{-\infty}^{i-1}}

\newcommand{\difi}{\delta_i(f)}

\newcommand{\difij}{\delta_{i+j}(f)}

\newcommand{\lip}{\textup{Lip}(E^{\Z},\R)}
\newcommand{\bx}{{\overline{x}}}
\newcommand{\by}{{\overline{y}}}
\newcommand{\bz}{{\overline{z}}}

\newtheorem{theorem}{{\small T}{\scriptsize HEOREM}}[section]
\newtheorem{corollary}{{\bf{\small C}{\scriptsize OROLLARY}}}[section]
\newtheorem{proposition}{{\bf{\small P}{\scriptsize ROPOSITION}}}[section]
\newtheorem{lemma}{{\bf{\small L}{\scriptsize EMMA}}}[section]
\newtheorem{remark}{{\bf{\small R}{\scriptsize EMARK}}}[section]
\newtheorem{definition}{{\bf{\small D}{\scriptsize EFINITION}}}[section]

\newcommand{\beq}{\begin{eqnarray}}
\newcommand{\eeq}{\end{eqnarray}}

\newcommand{\be}{\begin{equation}}
\newcommand{\ee}{\end{equation}}

\newcommand{\bl}{\begin{lemma}}
\newcommand{\el}{\end{lemma}}

\newcommand{\br}{\begin{remark}}
\newcommand{\er}{\end{remark}}

\newcommand{\bt}{\begin{theorem}}
\newcommand{\et}{\end{theorem}}

\newcommand{\bd}{\begin{definition}}
\newcommand{\ed}{\end{definition}}

\newcommand{\bp}{\begin{proposition}}
\newcommand{\ep}{\end{proposition}}

\newcommand{\bc}{\begin{corollary}}
\newcommand{\ec}{\end{corollary}}

\newcommand{\bpr}{\begin{proof}}
\newcommand{\epr}{\end{proof}}

\newcommand{\bi}{\begin{itemize}}
\newcommand{\ei}{\end{itemize}}

\newcommand{\ben}{\begin{enumerate}}
\newcommand{\een}{\end{enumerate}}

\def\1{{\mathchoice {\rm 1\mskip-4mu l} {\rm 1\mskip-4mu l}
{\rm 1\mskip-4.5mu l} {\rm 1\mskip-5mu l}}}

\begin{document}

\title{Concentration inequalities for Markov processes via coupling}
\author{Jean-Ren\'{e} Chazottes$^{\textup{{\tiny(a)}}}$, Frank Redig$^{\textup{{\tiny(b)}}}$\\
{\small $^{\textup{(a)}}$ Centre de Physique Th\'eorique, CNRS, \'Ecole Polytechnique}\\
{\small 91128 Palaiseau, France}\\
{\small $^{\textup{(b)}}$ Mathematisch Instituut Universiteit Leiden}\\{\small Niels Bohrweg 1, 2333 CA Leiden, The Netherlands}} 

\date{oct. 1, 2008}

\maketitle

\begin{abstract}
We obtain moment and Gaussian bounds for general coordinate-wise Lipschitz functions
evaluated along the sample path of a Markov chain.
We treat Markov chains on general (possibly unbounded) state spaces via a coupling method.
If the first moment of the coupling time exists, then we obtain
a variance inequality. If a moment of order $1+\eps$ of the
coupling time exists, then depending on the behavior of the stationary
distribution, we obtain higher moment bounds. This immediately implies
polynomial concentration inequalities.
In the case that a moment of order $1+\eps$ is finite uniformly in the
starting point of the coupling, we obtain a Gaussian bound.
We illustrate the general results with house of cards processes,
in which both uniform and non-uniform behavior of moments of the coupling time can occur.

\bigskip

\noindent{\bf Keywords}: Gaussian bound, moment bounds, house of cards process, Hamming distance.

\end{abstract}


\section{Introduction}

In this paper we consider a stationary Markov chain $X_n, n\in\Z$,
and want to obtain  inequalities for the probability that a function
$f(X_1,\ldots,X_n)$ deviates from its expectation.
In the spirit of concentration inequalities, one can try to bound
the exponential moment of $f-\E(f)$ in terms of the sum of squares of the
Lipschitz constants of $f$, as can be done in the case of independent random variables by several methods \cite{ledoux}. 

In the present paper, we want to continue the line of thought
developed in \cite{cckr,collet} where concentration inequalities
are obtained via a combination of martingale difference approach
(telescoping $f-\E(f)$) and coupling
of conditional distributions. In the case of an unbounded
state space, we cannot expect to find a coupling of
which the tail of the distribution of the coupling time can be controlled uniformly
in the starting points. This non-uniform dependence is thus
rather the rule than the exception and has to be dealt with if
one wants to go beyond the finite (or compact) state space situation.
Moreover, if the state space is continuous, then in general
two copies of the process cannot be coupled such that
they eventually coincide: we expect rather that in a coupling
the distance between the two copies can be controlled and
becomes small when we go further in time. We show that
a control of the distance suffices to obtain concentration inequalities.
This leads to a ``generalized coupling time'' which in discrete settings
coincides with the ordinary coupling time (in the case of a successful coupling).

In order to situate our results in the existing literature,
we want to stress that the main message of this paper
is the connection between the behavior of the generalized
coupling time and concentration inequalities. In order to illustrate
the possibly non-uniform behavior of the coupling time,
we concentrate on the simplest possible example of ``house
of cards'' processes (Markov chains on the natural numbers).
In this paper we restrict to the Gaussian concentration inequality
and moment inequalities. In principle, moment inequalities
with controll on the constants can be ``summarized'' in the
form of Orlicz-norm inequalities, but we do not want to
deal with this here.

The case of Markov chains  was first considered by Marton \cite{marton-1,marton0} : for uniformly contracting Markov chains, in particular
for ergodic Markov chains with finite state space, Gaussian concentration inequalities
are obtained.
The method developed in that paper is based on transportation cost-information inequalities. 
With the same technique, more general processes were considered by her in \cite{marton1}.
Later, Samson \cite{samson} obtained Gaussian concentration inequalities for some classes of Markov chains and
$\Phi$-mixing processes, by following Marton's approach.
Let us also mention the work by Djellout {\em et al.} \cite{arnaud} for further results in that direction.
Chatterjee \cite{chatterjee} introduced a version of Stein's method of exchangeable pairs to prove
Gaussian as well as moment concentration inequalities. Notice that moment inequalities were obtained
for Lipschitz functions of independent random variables in \cite{bouch}.
Using martingale differences,  Gaussian concentration inequalities were obtained in
\cite{k1,rio} for some classes of mixing processes. Markov contraction was used in \cite{k2} for ``Markov-type" processes ({\em e.g.}. hidden Markov chains).  

Related work to ours is  found in \cite{douc,douc2,douc3} where deviation or concentration inequalities \cite{douc} and speed of convergence
to the stationary measure \cite{douc2,douc3}  are obtained for subgeometric Markov chains, using a technique
of regeneration times and Lyapounov functions. 
Concentration properties
of suprema of additive functionals of Markov chains are
studied in \cite{adam}, using a technique of regeneration times.The example of the house of
cards process, and in particular its speed of relaxation to the stationary measure
is studied in \cite{douc2}, section 3.1. The speed of relaxation
to the stationary measure is of course related to the coupling time,
see {\em e.g.}. \cite{mao} for a nice recent account.
In fact, using an explicit coupling, we obtain concentration
inequalities in the different regimes of relaxation studied
in \cite{douc2}.  

Our paper is organized as follows. We start by defining the context
and introduce the telescoping procedure, combined with coupling.
Here the notion of coupling matrix is introduced. In terms
of this matrix we can (pointwise) bound the individual terms in the
telescopic sum for $f-\E(f)$. We then turn to the Markov case,
where there is a further simplification in the coupling matrix
due to the Markov property of the coupling.
In Section \ref{vario} we prove a variance bound
under the assumption that the first moment of the (generalized)
coupling time exists. In section
\ref{momo} we turn to moment inequalities. In this case
we require that a moment of order $1+\eps$ of the (generalized)
coupling time exists. This moment $M_{x,y, 1+\epsilon}$ depends on the starting
point of the coupling. The moment inequality for moments
of order $2p$ will then be valid if (roughly speaking)
the $2p$-th moment of $M_{x,y, 1+\epsilon}$ exists.
In Section \ref{gausio} we prove that if a moment of order $1+\eps$
of the coupling is finite, uniformly in the starting point, then
we have a Gaussian concentration bound.

Finally, Section \ref{exo} contains examples. In particular, we illustrate our approach
in the context of so-called house of cards processes, in which both
the situation of uniform case (Gaussian bound), as well
as the non-uniform case (all moments or moments up to a certain order)
are met. We end with application of our moment bounds to measure concentration
of Hamming neighborhoods and get non-Gaussian measure concentration bounds.

\section{Setting}

\subsection{The process}\label{tp}

The state space of our process is denoted by
$E$. It is supposed to be a metric space with distance $d$.
Elements of $E$ are denoted by $x,y,z$. $E$ is going
to serve as state space of a double sided
stationary process. Realizations of this process
are thus
elements of $E^{\Z}$ and are denoted by $\bx,\by,\bz$.

We denote by $(X_n)_{n\in\Z}$ a (two-sided)
stationary
process with values in $E$. The joint distribution
of $(X_n)_{n\in\Z}$ is denoted by $\pee$, and $\E$ denotes corresponding expectation.

$\Finf$ denotes the sigma-fields generated by $\{X_k: k\leq i\}$,
\[
\F_{-\infty}:=\bigcap_i \F_{-\infty}^i 
\]
denotes the tail sigma-field, and
\[
\F=\sigma\left(\bigcup_{i=-\infty}^\infty \F_{-\infty}^i\right). 
\]
We assume in the whole of this paper that $\pee$ is tail trivial, {\em i.e.},
for all sets $A\in \F_{-\infty}$, $\pee(A) \in \{ 0,1\}$.

For $i<j, i,j\in\Z$, we denote by $X_i^j$ the vector $(X_i, X_{i+1},\ldots, X_j)$, and
similarly we have the notation $X_{-\infty}^i$, $X_i^\infty$.
Elements of $E^{\{i,i+1,\ldots, j\}}$ ({\em i.e.}, realizations of $X_i^j$) are
denoted by $x_i^j$, and similarly we have $x_{-\infty}^i$, $x_i^\infty$.

\subsection{Conditional distributions, Lipschitz functions}\label{cl}

We denote by $\pee_{\xinf}$ the joint distribution of $\{X_j: j\geq i+1\}$ given $\Xinf=\xinf$.
We assume that this object is defined for all $\xinf$, {\em i.e.}, that there
exists a specification with which $\pee$ is consistent. This is automatically satisfied
in our setting, see \cite{goldstein}.

Further, $\co_{\xinf,\yinf}$ denotes a coupling of $\pee_{\xinf}$ and $\pee_{\yinf}$.


For $f:E^\Z\to \R$, we define the $i$-th Lipschitz constant

\[
\difi:=\sup\left\{\frac{f(\bx)-f(\by)}{d(x_i,y_i)} : x_j=y_j,\ \forall j\neq i,\ x_i\neq y_i\right\}.
\]

The function $f$ is said to be Lipschitz in the $i$-th coordinate if $\difi <\infty$,
and Lipschitz in all coordinates if $\difi <\infty$ for all $i$. We use
the notation $\delta(f)=(\delta_i(f))_{i\in\Z}$. We denote by
$\lip$ the set of all real-valued functions on $E^\Z$ which are Lipschitz
in all coordinates.

\section{Telescoping and the coupling matrix}

We start with $f\in\lip\cap L^1 (\pee)$, and
begin with the classical telescoping (martingale-difference) identity
$$
f-\E(f)= \sum_{i=-\infty}^\infty \Delta_i
$$
where
$$
\Delta_i:= \E(f\big| \Finf)-\E(f\big| \Finfi).
$$

We then write, using the notation of Section \ref{tp},
$$
\Delta_i = \Delta_i(\Xinf)= \int d\pee_{\Xinfi}(z_i)\times
$$
\begin{equation}\label{ding}
\int d\co_{\Xinf,\Xinfi z_i}(y_{i+1}^\infty,z_{i+1}^\infty)
\left[
f(\Xinf y_{i+1}^\infty)-f(\Xinfi z_{i}^\infty)
\right].
\end{equation}
For $f\in\lip$,  we have the following obvious telescopic
inequality
\begin{equation}\label{dong}
| f(\bx)- f(\by)| \leq \sum_{i\in\Z} \difi d(x_i,y_i).
\end{equation}
Combining \eqref{ding} and \eqref{dong} one obtains 
\begin{equation}\label{esope}
|\Delta_i(\Xinf)|\leq \sum_{j=0}^\infty D^{\Xinf}_{i,i+j} \difij
\end{equation}
where
\begin{equation}\label{coupmat}
D^{\Xinf}_{i,i+j}:= \int d\pee_{\Xinfi}(z_i)\int d\co_{\Xinf,\Xinfi z_i}(y_{i+1}^\infty,z_{i+1}^\infty) \thinspace d(y_{i+j},z_{i+j}).
\end{equation}
This is an upper-triangular random matrix which we call the coupling matrix associated with the process $(X_n)$ and $\co$,
the coupling of the conditional distributions. 
As we obtained before in \cite{cckr}, in the context of $E$ a finite set,
the decay properties of the matrix elements $D^{\Xinf}_{i,i+j}$
({\em i.e.}, how these matrix elements become small when $j$ becomes large) determine
the concentration properties of $f\in\lip$, via the control
\eqref{esope} on $\Delta_i$, together with
Burkholder's inequality \cite[Theorem 3.1, p. 87]{burk1}, which relates
the moments of $f-\E(f)$ with powers of the sum of
squares of $\Delta_i$.
The non-uniformity (as a function of  the realization of $\Xinf$) of the decay of the matrix elements
as a function of $j$ (which we encountered {\em e.g.}. in the low-temperature Ising model \cite{cckr}) will be {\em typical}
as soon as the state space $E$ is unbounded. Indeed, if
starting points in the coupling are further away, then
it takes more time to get the copies close in the coupling .

\br
The same telescoping procedure can be
obtained for ``coordinate-wise H\"{o}lder''
functions, {\em i.e.}, functions such that 
for some $0< \alpha <1$
\[
\delta^\alpha_i (f):=\sup\left\{\frac{f(\bx)-f(\by)}{d^\alpha(x_i,y_i)} : x_j=y_j,\ \forall j\neq i,\ x_i\neq y_i\right\}.
\]
is
finite for all $i$. In \eqref{coupmat}, we
then have to replace $d$ by
$d^\alpha$.
\er

\section{The Markov case}

We now consider $(X_n)_{n\in\Z}$ to be a stationary and ergodic Markov chain. We denote by 
$p(x,dy):=\pee(X_1\in dy\big| X_0=x)$ the transition kernel. 
We let $\nu$ be the unique stationary measure of the Markov chain.
We denote by $\pee_\nu$ the path space measure
of the stationary process $(X_n)_{n\in\Z}$.
By $\pee_x$ we denote the distribution of $(X_1^\infty)$, for
the Markov process conditioned on $X_0=x$.
 
We further suppose that the coupling $\co$ of Section \ref{cl} is Markovian, and denote
by $\co_{x,y}$ the coupling started from $x,y$, and corresponding
expectation by $\hat{\E}_{x,y}$. More precisely, by the Markov
property of the coupling we then have that 
$$
\co_{\xinf,\yinf}=\co_{x_i,y_i}
$$
is a Markovian coupling $((\hat{X}^{(1)}_n,\hat{X}^{(2)}_n))_{n\in\N}$ of the Markov chains $(X_n)_{n\geq 0}$ starting from $X_0=x_i$, resp. $Y_0=y_i$.
In this case the expression \eqref{coupmat} of the coupling matrix simplifies to
\[
\Psi_{X_{i-1}, X_i} (j):=D^{X_{i-1},X_i}_{i,i+j}=
\int p(X_{i-1},dy) \int d\co_{X_i,y}(u_0^\infty,v_0^\infty) \thinspace d(u_{j},v_{j})
\]
With this notation, \eqref{esope} reads
\be\label{esopsi}
|\Delta_i (X_{i-1}, X_{i})| \leq \sum_{j\geq 0} \Psi_{X_{i-1}, X_i} (j) \delta_{i+j} f.
\ee

We define the ``generalized coupling time''
\begin{equation}\label{defT}
\T(u_0^\infty,v_0^\infty):= \sum_{j=0}^\infty d(u_{j},v_{j}).
\end{equation}
In the case $E$ is a discrete (finite or countable) alphabet, the ``classical'' coupling time
is defined as usual
\[
T (u_0^\infty,v_0^\infty): = \inf\{ k\geq 0: \forall j\geq k : u_j=v_j\}.
\]
If we use the trivial distance $d(x,y)=1$ if $x\neq y$ and $d(x,y)=0$ if $x=y$, for $x,y\in E$, then
we have
\be\label{pipo}
 d(u_{j},v_{j})\leq \1\{T(u_0^\infty,v_0^\infty) \geq j\}
\ee
and hence
\[
\T(u_0^\infty,v_0^\infty) \leq T(u_0^\infty,v_0^\infty).
\]
Of course, the same inequality remains true if $E$ is a bounded
metric space with $d(x,y)\leq 1$ for $x,y\in E$. However
a ``successful coupling'' ({\em i.e.}, a coupling with $T<\infty$)
is not expected to exist in general in the case
of a non-discrete state space. It can however exist, see
{\em e.g.}.\ \cite{fey} for a successful coupling in the
context of Zhang's model of self-organized criticality.
Let us also mention that the ``generalized coupling time" unavoidably appears in the context of dynamical systems \cite{collet}.

In the discrete case, using \eqref{esopsi} and \eqref{pipo}, we
obtain  the following inequality: 
\be\label{guerrier}
\Psi_{x,y} (j) \leq  \sum_{z} p(x,z) \hat{\E}_{z,y} \left(\1\{T(u_0^\infty,v_0^\infty) \geq j\}\right)
\ee
whereas in the general (not necessarily discrete) case we have, by \eqref{defT}, and monotone convergence,
\be\label{zann}
\sum_{j\geq 0}\Psi_{x,y} (j) = \int p(x, dz) \hat{\E}_{y, z} (\T)
\leq \int p(x, dz) \hat{\E}_{y, z} (T).
\ee

\br

So far, we made a telescoping of $f-\E(f)$ using an increasing family of sigma-fields. One can 
as well consider a decreasing family of sigma-fields, such as $\F_i^\infty$, defined to be the
sigma-fields generated by $\{X_k: k\geq i\}$. 
We then have,  {\em mutatis mutandis}, the same inequalities
using ``backward telescoping''
$$
f-\E(f)=\sum_{i=-\infty}^\infty \Delta_i^*,
$$ 
where
$$
\Delta_i^*:=\E(f\big|\F_i^\infty)-\E(f\big|\F_{i+1}^\infty).
$$
and estimating $\Delta_i^*$ in a completely parallel way, by introducing a lower-triangular
analogue of the coupling matrix
matrix.

Backward telescoping is natural in the context of dynamical systems where the forward 
process is deterministic, hence cannot be coupled (as defined above) with two different initial conditions
such that the copies become closer and closer.
However, backwards in time, such processes are non-trivial Markov chains for which a coupling can
be possible with good decay properties of the coupling matrix. See \cite{collet} for a concrete example
with piecewise expanding maps of the interval.
\er

\section{Variance inequality}\label{vario}

For a real-valued sequence $(a_i)_{i\in\Z}$, we denote the usual $\ell_p$-norm by
\[
 \| a\|_p = \left(\sum_{i\in\Z} |a_i|^p \right)^{1/p}.
\]
Our first result concerns the variance of a $f\in\lip$.
\begin{theorem}\label{DevroyeMarkovchains}
Let $f\in \lip\cap L^2 (\pee_\nu)$. Then 
\begin{equation}\label{devroye}
\var(f)\leq 
C\|\delta(f)\|_2^2
\end{equation}
where
$$
C =  \int \nu (dx) \times
$$
\be\label{devconst}
\int p(x, dz) \int p(x, dy)\int p(x, du) \thinspace \hat{\E}_{z,y} (\T) \hat{\E}_{z,u} (\T).
\ee
As a consequence, we have the concentration inequality
\be\label{convar}
\forall t>0, \quad \pee(|f-\E(f)|\geq t)\leq C\ \frac{\|\delta(f)\|_2^2}{t^2}.
\ee
\end{theorem}

\bigskip

\begin{proof}
We estimate, using \eqref{esopsi} and stationarity
\[
\E(\Delta_i^2 )\leq
\E\left(\sum_{j\geq 0} \Psi_{X_0,X_1} (j) \delta_{i+j} (f)\right)^2=  \E\left(\left(\Psi_{X_0, X_1} *\delta (f)\right)_i^2\right).
\]
where $*$ denotes convolution, and where we extended $\Psi$ to $\Z$
by putting it equal to zero for negative integers.
Since 
$$
\var(f)=\sum_{i=-\infty}^\infty \E\left(\big(\Delta_i\big)^2\right)
$$
Using Young's inequality, we then obtain,
$$
\var(f) \leq \E\left(\big\|\Psi_{X_0,X_1}*\delta(f)\big\|_2^2\right) \leq \E\left(\big\|\Psi_{X_0,X_1}\big\|_1^2\right)\thinspace \big\|\delta(f)\big\|_2^2.
$$
Now, using the equality in \eqref{zann}
\beq
&&\E\left(\big\|\Psi_{X_0,X_1}\big\|_1^2\right)
\nonumber\\
&=&\E \left(\int p(X_0, dy) \hat{\E}_{X_1, y} (\T)\right)^2
\nonumber\\
&=&
\int \nu (dx) p(x, dz) \left( \int p(x, dy) \hat{\E}_{z,y} (\T)\right)^2
\nonumber\\
&= &
\nonumber
\int \nu (dx) \int p(x, dz) \int p(x, dy)\int p(x, du)\thinspace \hat{\E}_{z,y} (\T) \hat{\E}_{z,u} (\T),
\eeq
which is \eqref{devroye}.
Inequality \eqref{convar} follows from Chebychev's inequality.
\end{proof}

The expectation in \eqref{devconst} can be interpreted as follows.
We start from a point $x$ drawn from the
stationary distribution and generate three independent
copies $y,u,z$ from the Markov chain at time
$t=1$ started from $x$. With these initial points we start the
coupling in couples $(y,z)$ and $(u,z)$, and compute the expected
coupling time.

\section{Moment inequalities}\label{momo}

In order to control higher moments of $(f-\E(f))$, we have
to tackle higher moments of the sum $\sum_{i}\Delta_i^2$
and for these we cannot use the simple
stationarity argument used in the estimation of the variance.

Instead, we start again from
\eqref{esopsi} and let $\lambda^2 (j) := (j+1)^{1+\epsilon}$ where
$\epsilon >0$.

We then obtain, using Cauchy-Schwarz inequality:
\beq
\nonumber
|\Delta_i| & \leq & \sum_{j\geq 0} \lambda (j) \Psi_{X_{i-1}, X_i} (j) \frac{\delta_{i+j} (f)}{\lambda (j)}
\nonumber\\
&\leq &
\left(\sum_{j\geq 0} \lambda (j)^2 \left(\Psi_{X_{i-1}, X_i} (j)\right)^2
\sum_{k\geq 0}\left(\frac{\delta_{i+k} (f)}{\lambda (k)}\right)^2\right)^{1/2}.
\nonumber
\eeq
Hence
\be\label{koko}
\Delta_i^2 \leq \Psi_\eps^2 (X_{i-1},X_i) \left((\delta(f))^2*\frac{1}{\lambda^2}\right)_i
\ee
where $\delta (f)^2$ denotes the sequence with components $\left(\delta_i (f)\right)^2$, and
where
\be\label{psipsi}
\Psi^2_\eps (X_{i-1},X_i)= \sum_{j\geq 0} (j+1)^{1+\eps} \left(\Psi_{X_{i-1}, X_i} (j)\right)^2.
\ee
Moment inequalities will now be expressed in terms of moments of $\Psi_\eps^2$.

\subsection{Moment inequalities in the discrete case}

We first deal with a discrete state space $E$. Recall \eqref{pipo}.

\bl
In the discrete case, {\em i.e.}, if $E$ is a countable
set with the discrete metric, then, for all $\eps>0$,  we have the estimate
\beq\label{psila}
&&\Psi_\eps^2 (X_{i-1}, X_i) 
\nonumber\\
&&\leq \frac12 \left(\sum_{z} p(X_{i-1}, z)  \hat{\E}_{X_i,z} ((T+1)^{1+\frac{\epsilon}{2}})\right)^2.
\eeq
\el

\bpr
Start with
\beq
\Psi_\eps^2 &=& \sum_{j\geq 0} (j+1)^{1+\eps}\left(\Psi_{X_{i-1}, X_i} (j)\right)^2
\nonumber\\
&\leq &
\sum_{z,u}\sum_{j\geq 0} (j+1)^{1+\eps} p(X_{i-1}, z) p(X_{i-1}, u)\hat{\pee}_{X_i,z}(T\geq j)
\hat{\pee}_{X_i,u}(T\geq j).
\nonumber
\eeq
Proceed now with
\beq
&&\sum_{j\geq 0} (j+1)^{1+\eps} p(X_{i-1}, z) p(X_{i-1}, u)\hat{\pee}_{X_i,z}(T\geq j)
\hat{\pee}_{X_i,u}(T\geq j)
\nonumber\\
&=& \sum_{k=0}^\infty\sum_{l=0}^\infty\sum_{j=0}^{l\wedge k} (j+1)^{1+\eps}\ \hat{\pee} (T_1=k, T_2=l)
\nonumber\\
& \leq & 
\frac12   
\sum_{k=0}^\infty\sum_{l=0}^\infty (l\wedge k+1)^{2+\eps}\ \hat{\pee} (T_1=k, T_2=l)
\nonumber\\
&=&
\nonumber
\hat{\E} \left(((T_1+1)\wedge( T_2+1))^{2+\eps}\right),
\eeq
where we denoted by $T_1$ and $T_2$ two independent
coupling times corresponding to two independent copies of the coupling started
from $(X_i,z)$, resp.\ $(X_i,u)$.

Now use that for two independent non-negative real-valued random variables we have
\[
\E\big( (X\wedge Y)^{2+\eps}\big)\leq \E(X^{1+\frac{\eps}{2}})\E(Y^{1+\frac{\eps}{2}}).
\]
The lemma is proved.
\epr
In order to arrive at moment estimates, we want an estimate
for
$\E\big(\sum_i \Delta_i^2\big)^p$. This is the content of the next lemma. We denote, as usual,
$\zeta (s) =\sum_{n=1}^{\infty}(1/n)^s$.
\bl
For all $\eps>0$ and integers $p>0$ we have 
\beq\label{mompou}
\E\big(\sum_i \Delta_i^2\big)^p
&\leq &
\left(\frac{\zeta (1+\eps)}{2}\right)^p \|\delta (f)\|_2^{2p} \sum_{x,y} \nu (x) p(x,y) \times
\nonumber\\
&&  \left(\sum_{z}p(x, z)  \hat{\E}_{y,z} ((T+1)^{1+\frac{\epsilon}{2}}))\right)^{2p}.
\eeq
\el
\bpr
We start from
\beq
&&\E(\sum_i \Delta_i^2)^p
\nonumber\\
&&
\nonumber
\leq \sum_{i_1,\ldots, i_p} \E\left(\prod_{l=1}^p\Psi_\eps (X_{i_l-1},X_{i_l})^2\right) \prod_{l=1}^p\left((\delta(f))^2*\frac{1}{\lambda^2}\right)_{i_l}.
\eeq
Then use H\"{o}lder's inequality and stationarity, to obtain
\beq
\E\big(\sum_i \Delta_i^2\big)^p
&\leq &
\E (\Psi_\eps^{2p} (X_0, X_1))\times \Big\|(\delta (f))^2*\frac{1}{\lambda^2}\Big\|_1^p
\nonumber\\
&\leq &
\E (\Psi_\eps^{2p} (X_0, X_1))\times\Big\|\frac{1}{\lambda^2}\Big\|_1^p \ \|(\delta (f))^2\|_1^p
\nonumber\\
&=&
\nonumber
\E (\Psi_\eps^{2p} (X_0, X_1))\times\Big\|\frac{1}{\lambda^2}\Big\|_1^p \ \|(\delta (f))\|_2^{2p}
\eeq
where in the second inequality we used Young's inequality.
The lemma now follows from 
\eqref{psila}.
\epr
We can now formulate our moment estimates in the discrete case.

\bt\label{HigherMomentsMarkovchains}
Suppose $E$ is a countable set with discrete metric.
Let $p\geq 1$ be an integer and $f\in\lip\cap L^{2p}(\pee)$.
Then for all $\eps >0$ we have the estimate
\be\label{kloot}
\E (f-\E(f))^{2p}  \leq C_{p}  \|\delta (f)\|_2^{2p}
\ee
where
\beq\label{momco}
C_{p} &= &
(2p-1)^{2p}\left(\frac{\zeta (1+\eps)}{2}\right)^p  \times
\nonumber\\
&&  \sum_{x,y} \nu (x) p(x,y) \left(\sum_{z}p(x, z) \hat{\E}_{y,z} \big((T+1)^{1+\frac{\epsilon}{2}}\big)\right)^{2p}.
\eeq
As a consequence we have the concentration inequalities
\be\label{polyci}
\forall t>0, \quad \pee(|f-\E(f)|\geq t)\leq C_p\ \frac{\|\delta(f)\|_2^{2p}}{t^{2p}}.
\ee
\et

\bigskip

\bpr
By Burkholder's inequality \cite[Theorem 3.1, p. 87]{burk1}, one gets
$$
\E\left(\big(f-\E(f)\big)^{2p}\right)
\leq 
(2p-1)^{2p} \thinspace \E\Big(\big(\sum_i \Delta_i^2\big)^p\Big)
$$
and \eqref{momco} then follows from \eqref{mompou}, whereas \eqref{polyci} follows from \eqref{momco} by Markov's inequality.
\epr
\br\label{strounf}
Theorem \ref{HigherMomentsMarkovchains} for $p=1$ is weaker than
Theorem \ref{DevroyeMarkovchains}: indeed, for \eqref{devroye}
to hold
we only need to have
the first moment of the coupling time to be finite. 
\er
\br
A typical behavior (see the examples below) of the coupling time is as follows:
\[
 \co_{x,y} (T\geq j) \leq C(x,y) \phi (j)
\]
Here $C(x,y)$ is a constant that depends, in general in an unbounded
way, on the starting points
$(x,y)$ in the coupling
and where $\phi (j)$,  determining the tail of the coupling
time does not depend on the starting points.
Therefore, for the finiteness of the constant $C_p$ in \eqref{momco}
we need that the tail-estimate $\phi (j)$ decays fast enough so that $\sum_{j} j^\eps \phi (j) <\infty$
(that does not depend on $p$), and next the $2p$-th power of the constant
$C(x,y)$ has to be integrable (this depends on $p$).
\er

\subsection{The general state space case}

In order to formulate the general state space version of these results,
we introduce the expectation 
\[
 \tilde{\E}_{x,y} ( F (\overline{u},\overline{v})) = \int p(x,dz) 
\int \hat{\pee}_{y,z} (d\overline{u}, d\overline{v}) F (\overline{u},\overline{v}).
\]
We can then rewrite
\[
 \Psi_\eps^2 (x,y) = \sum_{j\geq 0} (j+1)^{1+\eps}\left(\tilde{\E}_{x,y} d(u_j, v_j)\right)^2.
\]
We introduce
\[
 \alpha^{x,y}_j =\left(\tilde{\E}_{x,y} (d(u_{j}, v_{j}))- \tilde{\E}_{x,y} (d(u_{j+1}, v_{j+1}))\right).
\]
This quantity is the analogue of $\hat{\pee} (T=j)$ of the discrete case.
We then define
\be\label{mr}
 M_r^{x,y} = \sum_{j\geq 0} (j+1)^r \alpha_j^{x,y}
\ee
which is the analogue of the $r$-th moment of the coupling time.
The analogue of Theorem \ref{HigherMomentsMarkovchains} then becomes the following.

\bt
Let $p\geq 1$ be an integer and $f\in\lip\cap L^{2p}(\pee)$.
Then for all $\eps >0$ we have the estimate
\[
\E (f-\E(f))^{2p}  \leq C_{p}  \|\delta (f)\|_2^{2p}
\]
where
\[
C_{p} = 
(2p-1)^{2p}\left(\frac{\zeta (1+\eps)}{2}\right)^p  
\int\nu (dx) p(x,dy) \left(M^{x,y}_{1+\frac{\eps}{2}}\right)^{2p}.
\]
\et

\section{Gaussian concentration bound}\label{gausio}

If one has a uniform estimate of the quantity
\eqref{psipsi}, we obtain a corresponding uniform estimate
for $\Delta_i^2$, and via Hoeffding's inequality, a
Gaussian bound for $ (f-\E(f))$. This is formulated in the
following theorem.
\bt
Let $E$ be a countable set with the discrete metric.
Let $f\in \lip$ such that $\exp (f)\in L^1 (\pee)$. Then for any $\eps >0$
we have
\be\label{gogo}
\E\left( e^{f-\E(f)}\right) \leq e^{C_\eps \|\delta (f)\|_2^2/16}
\ee
where
\be\label{gogoce}
C_\eps = \zeta (1+\eps) \left(\sup_{u,v} \hat{\E}_{u,v} (T^{1+\eps})\right)^2.
\ee
The general state space analogue of this bound is obtained by replacing
$\hat{\E}_{u,v} (T^{1+\eps})$ by $M_{1+\eps}^{u,v}$ in \eqref{gogoce} (where $M_r^{x,y}$ is defined in \eqref{mr}).
\et 
\bpr
From \eqref{psila} we get
\beq
\Psi_{\eps}^2 (x,y) &\leq &
\frac12 \sum_{z,u} p(x,z) p(x,u) \hat{\E}_{y,z} (T^{1+\frac{\eps}{2}}) \hat{\E}_{y,u} (T^{1+\frac{\eps}{2}})
\nonumber\\
&\leq &
\frac12 \left(\sup_{u,v}\hat{\E}_{u,v} (T^{1+\frac{\eps}{2}})\right)^2.
\nonumber
\eeq
Now start from the classical Azuma-Hoeffding inequality \cite{mac}
\[
\E\left( e^{f-\E(f)}\right)\leq e^{\frac18 \sum_i \|\Delta_i\|^2_\infty}.
\]
Therefore, we estimate, using \eqref{koko} and Young's inequality
\[
\sum_{i} \|\Delta_i\|^2_\infty \leq \frac12 \zeta (1+\eps) \left(\sup_{u,v}\hat{\E}_{u,v} (T^{1+\frac{\eps}{2}})\right)^2
\|\delta (f)\|_2^2
\]
which establishes \eqref{gogo}. 
\epr

\br
Let us take $f$ bounded Lipschitz and depending on finitely many coordinates, that is, assume
there exists $n\geq 1$ such that $f(\bx)=f(\by)$ whenever $x_i=y_i$, $i=-n,\ldots,n$.
Since $\E(e^{\lambda f})<\infty$ for every $\lambda>0$, inequality \eqref{gogo} implies a concentration inequality by the optimized exponential Chebychev inequality:
\be\label{concegauss}
\forall t> 0, \quad \pee(|f-\E(f)|\geq t)\leq 2\exp\left(\frac{-4C_\eps t^2} {\|\delta(f)\|_2^2} \right)\cdot
\ee
More generally  we may have $\E(e^{\lambda f})<\infty$ for $\lambda\in(0,\lambda_0]$, for some $\lambda_0>0$.
Then there are two regimes: for $t\leq t_0$, the previous inequality holds, but for $t>t_0$, the bound is
of the form $e^{-ct}$.
\er

\br
The assumption that a moment of order $1+\epsilon$ of the coupling time exists,
which is uniformly bounded in the starting point,
can be weakened to the same property for the first moment, if we have
some form of monotonicity.
More precisely, we say that a coupling has the monotonicity property, if
there exist ``worse case starting points'' $x_u, x_l$, which have
the property that
\[
\sup_{x,y}\hat{\pee}_{x,y} (T\geq j) \leq \hat{\pee}_{x_u, x_l} (T\geq j)
\]
for all $j\geq 0$.
In that case, using \eqref{guerrier}, we can start from \eqref{esopsi}
and obtain, in the discrete case, the uniform bound
\[
\Delta_i \leq \sum_{j\geq 0} \hat{\pee}_{x_u, x_l} (T\geq j) \delta_{i+j} f
\]
and via Azuma-Hoeffding inequality, combined with Young's inequality, we then obtain the Gaussian
bound
\eqref{gogo} with
\[
 C= \frac12 \hat{\E}_{x_u, x_l} (T).
\]
Finally, it can happen (especially if the state space is unbounded)
that the coupling has no worst case starting
points, but there is a sequence $x^n_u, x^n_l$ of elements
of the state space such that $\hat{\pee}_{x^n_u, x^n_l} (T\geq j)$
is a non-decreasing sequence in $n$ for every fixed $j$ and
\[
 \hat{\pee}_{x,y} (T\geq j) \leq \lim_{n\to\infty}\hat{\pee}_{x^n_u, x^n_l} (T\geq j).
\]
\textup{(}{\em E.g.}.,  in the case of the state space $\Z$, we can think of
the sequence $x^n_u\to \infty$ and $x_l^n\to -\infty$.\textup{)}
In that case, from monotone convergence we have the Gaussian
concentration bound with
\[
 C= \lim_{n\to\infty}\frac12 \hat{\E}_{x^n_u, x^n_l} (T).
\]
\er

\section{Examples}\label{exo}

\subsection{Finite-state Markov chains}

As we mentioned in the introduction, this case was already considered by K. Marton (and others), but it illustrates our method in the most simple setting, and gives also an alternative proof
in this setting.

Indeed, if the chain is aperiodic and irreducible, then it is well-known
\cite{iceman}, 
$$
\sup_{u,v\in E} \hat{\pee}_{u,v}(T\geq j)\leq c\ \rho^j
$$
for all $j\geq 1$ and some $c>0$.
Hence  the Gaussian bound \eqref{gogo} holds.

\subsection{House of cards processes}\label{hoc}

These are Markov chains on the set of natural numbers which
are useful in the construction of couplings for processes
with long-range memory, and dynamical systems, see {\em e.g.}.\
\cite{bfg}.

More precisely, a house of cards process is a Markov chain on the natural numbers
with transition probabilities
\[
\pee (X_{k+1}= n+1|X_k =n) = 1-q_n= 1- \pee (X_{k+1}=0|X_k=n),
\]
for $n=0,1,2,\ldots$,
{\em i.e.}, the chain can go ``up" with one unit or go ``down" to zero.
Here, $0<q_n<1$.

In the present paper, house of card chains serve as a nice class of
examples where we can have moment inequalities up to a certain
order, depending on the decay of $q_n$, and even Gaussian
inequalities. Given a sequence of independent uniformly
distributed random variables $(U_k)$ on $[0,1]$, we can view
the process $X_k$ generated via the recursion
\be\label{rec}
X_{k+1}= (X_k +1) \1\{U_{k+1}\geq q_{X_k}\}.
\ee
This representation also yields a coupling of the process
for different initial conditions. The coupling has the property
that when the coupled chains meet, they stay together forever.
In particular, they will stay together forever after they hit
together zero.
For this coupling, we have the following estimate.
\bl
Consider the coupling defined via \eqref{rec}, started
from initial condition $(k,m)$ with $k\geq m$. Then we have
\be\label{coupes}
\hat{\pee}_{k,m} (T\geq t) \leq \prod_{j=0}^{t-1} (1-q^*_{k+j})
\ee
where
\[
q^*_n = \inf_{s\leq n} q_s.
\]
\el
\bpr
Call $Y^k_t$ the process defined by \eqref{rec} started from
$k$, and define $Z^k_t$, a process started from $k$ defined
via the recursion
\[
Z_{t+1} = (Z_t+1) \1\left\{U_{t+1} \geq q^*_{Z_t}\right\},
\]
where $U_t$ is the same sequence of independent
uniformly distributed random variables as in \eqref{rec}.
We claim that, for all $t\geq 0$,
\[
Y^k_t\leq Z^k_t, Y^m_t \leq Z^k_t.
\]
Indeed, the inequalities hold at time zero. Suppose they hold at
time $t$, then, since $q^*_n$ is non-increasing as a function of $n$,
\[
q_{Y^k_t} \geq q^*_{Y^k_t}\geq q^*_{Z^k_t} \quad \textup{and} \quad q_{Y^m_t}\geq q^*_{Y^m_t}\geq q^*_{Z^k_t}
\]
whence
\[
 \1\left\{U_{t+1}\geq q^*_{Z^k_t}\right\}\geq \1\left\{U_{t+1}\geq q_{Y^k_t}\right\}
\ \mbox{and} \ \1\left\{U_{t+1}\geq  q^*_{Z^k_t}\right\} \geq \1\left\{U_{t+1}\geq q_{Y^m_t}\right\}.
\]
Therefore, in this coupling, if $Z^k_t=0$, then  $Y^m_t=Y^k_t=0$, and hence
the coupling time is dominated by the first visit of $Z^k_t$ to zero, which gives
\[
 \hat{\pee}_{k,m} (T\geq t)\leq \pee (Z^k_n\not=0, n=1,\ldots, t-1)= \prod_{j=1}^{t-1} (1-q^*_{k+j}).
\]
\epr
The behavior \eqref{coupes} of the coupling time shows the typical non-uniformity
as a function of the initial condition. More precisely,
the estimate in the rhs of \eqref{coupes} becomes bad for large
$k$.
We now look at three more concrete cases.

\begin{enumerate}

\item {\bf Case 1:}
\[
q_n = \frac{1}{n^\alpha} \ , \; n\geq 2, \; 0< \alpha <1.
\]
Then it is easy to deduce from \eqref{coupes} that
\be\label{alca}
\hat{\pee}_{k,m} (T\geq t)\leq C \exp \left(-\frac{1}{1-\alpha} \left( (t+k)^{1-\alpha}-k^{1-\alpha}\right)\right).
\ee
The stationary (probability) measure is given by:
\be\label{statmeas}
\nu ( k) = \pi_0 c_k
\ee
with
\be\label{statbul}
 c_k = \prod_{j=0}^k (1-q_j)
\ee
which is bounded from above by
\be\label{vende}
 c_k \leq C' \exp \left(-\frac{1}{1-\alpha} k^{1-\alpha}\right).
\ee

>From \eqref{alca}, combined with \eqref{statmeas}, \eqref{vende}, it is then easy to see
that the constant $C_p$ of \eqref{momco} is finite for all $p\in \N$.
Therefore, in that case the moment inequalities 
\eqref{kloot} hold, for all $p\geq 1$.

\item {\bf Case 2:}
\[
q_n = \frac{\gamma}{n}\quad (\gamma>0)
\]
for $n\geq \gamma+1$, and other values $q_i$ are arbitrary.
In this case we obtain from \eqref{coupes} the estimate
\[
\hat{\pee}_{k,m} (T\geq t)\leq C_\gamma\frac{(k+1)^\gamma}{(k+t)^\gamma}
\]
and for the stationary measure we have \eqref{statbul} with
\[
c_k\leq C'_\gamma k^{-\gamma}.
\]
The constant $C_p$ of \eqref{momco}
is therefore bounded by
\[
 C_p \leq C_\gamma C'_\gamma C_1 C_2
\]
where $C_1= (2p-1)^{2p} (\zeta(1+\epsilon)/2)^p$ is finite independent of $\gamma$, and where
\[
 C_2= C_2 (p)
\leq \sum_{k\geq 1} k^{-\gamma} \left(\hat{\E}_{k+1,0} (T+1)^{1+\frac{\epsilon}{2}}\right)^{2p}
\]
so we estimate
\[
\hat{\E}_{k+1,0} (T+1)^{1+\frac{\epsilon}{2}}
\leq (1+\delta)\ \sum_{t=0}^\infty \frac{(k+1)^\gamma (t+1)^\delta}{(t+k)^\gamma},
\]
where $\delta:=\epsilon/2$.
To see when $C_2<\infty$, we first look at the behavior of
\[
I(a,b,k):= \sum_{t=1}^\infty \frac{t^a}{(t+k)^b} .
\]
The sum in the rhs is convergent for $ b-a>1$, in which case it
behaves as $k^{1+ a-b}$ for $k$ large, which gives
for our case $a= \delta$, $b=\gamma$, $\gamma >1+\delta$. In that case, we
find that $C_2(p)$ is finite as soon as
\[
 \sum_{k\geq 1} k^{-\gamma + 2p\gamma} k^{2p(\delta -\gamma+1)} <\infty
\]
which gives
\[
 \gamma > 1+ 2p \delta+2p.
\]
Hence, in this case, for $\gamma >1+\delta$, we obtain the moment estimates \eqref{HigherMomentsMarkovchains} up to order
$p < (\gamma-1)/2(\delta+1)$.

\item {\bf Case 3:}
\[
q:= \inf \{ q_n:n\in\N\} >0 
\]
then
we have the uniform estimate
\[
\sup_{k,m}\hat{\pee}_{k,m} (T\geq t) \leq (1-q)^t
\]
which gives the Gaussian
concentration bound
\eqref{gogo}
with
$C= \frac{1}{2(1-q)}$.
\end{enumerate}
\subsection{Ergodic interacting particle systems}
As a final example, we consider spin-flip dynamics in the so-called
$M<\epsilon$ regime. These are Markov processes on the space
$E= \{0,1\}^{S}$, with $S$ a countable set.
This is a metric space with distance
\[
 d(\eta, \xi) = \sum_{n=1}^\infty 2^{-n} |\eta_{i_n}-\xi_{i_n}|
\]
where $n\mapsto i_n$ is a bijection from $\N$ to $S$.

The space $E$ is interpreted as set of configurations of
``spins'' $\eta_i$ which can be up ($1$) or down $(0)$ and are defined
on the set $S$ (usually taken to be a lattice such as $\Z^d$). The
spin at site $i\in S$ flips at a configuration dependent rate $c(i,\eta)$.
The process is then defined via its generator on local functions defined
by
\[
 L f(\eta) = \sum_{i\in S} c(i,\eta) (f(\eta^i) -f(\eta))
\]
where $\eta^i$ is the configuration $\eta$ obtained from
$\eta$ by flipping at site $i$. See \cite{ligg} for
more details about existence and ergodicity of such processes.

We assume here that we are in the so-called ``$M<\epsilon$ regime", where
we have the existence of a coupling (the so-called ``basic coupling") for which
we have the estimate
\be\label{basiccoup}
 \hat{\pee}_{\eta^j, \eta}( \eta_i (t) \not= \zeta_i (t))< e^{-\epsilon t} e^{\Gamma (i,j) t}
\ee
with $\Gamma (i,j)$ a matrix indexed by $S$ with finite $\ell_1$-norm $M<\epsilon$.
As a consequence, from any initial configuration, the system evolves
exponentially fast to its unique equilibrium measure which we denote $\mu$.
The stationary Markov chain is then defined as
$X_n = \eta_{n\delta}$ where $\delta>0$, and
$\eta_0= X_0$
is distributed according to $\mu$.

In the basic coupling, from  \eqref{basiccoup}, we obtain the uniform estimate
\[
 \tilde{\E}_{\eta,\zeta} (d(\eta (k), \zeta(k))\leq e^{-{(M-\epsilon)}k\delta}.
\]
As a consequence, the quantity $M_r^{\eta,\zeta}$ of \eqref{mr} is finite
uniformly in $\eta,\zeta$, for every $r>0$.
Therefore, we have the Gaussian bound \eqref{gogo}
with 
\[
 C \leq \sum_{k\geq 1}k^{1+\epsilon} e^{-(M-\epsilon)k\delta} <\infty.
\]

\subsection{Measure concentration of Hamming neighborhoods}

We apply Theorem \ref{HigherMomentsMarkovchains} to measure concentration of Hamming neighborhoods. The case of
contracting Markov chains was already (and first) obtained in \cite{marton-1} as a consequence of an information divergence
inequality.  We can easily obtain such Gaussian measure concentration from \eqref{gogo}. But, by a well-known result of
Bobkov and G\"otze \cite{bg}, \eqref{gogo} and that information divergence inequality are in fact equivalent.  The interesting situation
is when \eqref{gogo} does not hold but only have moment bounds.  

Let $A,B\subset E^n$ be two sets and denote by $\bar{d}(A,B)$ their normalized Hamming distance,
$\bar{d}(A,B)=\inf\{\bar{d}(x_1^n,y_1^n):x_1^n\in A,y_1^n\in B\}$,
where 
$$
\bar{d}(x_1^n,y_1^n)=\frac{1}{n}\sum_{i=1}^n d(x_i,y_i),
$$
$d(x_i,y_i)=1$ if $x_i\neq y_i$, and $0$ otherwise. The $\varepsilon$-neighborhood of $A$ is then
$$
[A]_\varepsilon=\{y_1^n: \inf_{x_1^n\in A} \bar{d}(x_1^n,y_1^n)\leq \varepsilon\}.
$$
\bt
Take any $n\in\N$ and let $A\subset E^n$ a measurable set with $\pee(A)>0$. Then, under the assumptions of
Theorem \ref{HigherMomentsMarkovchains}, we have, for all $p\geq 1$,
\[
\pee([A]_\varepsilon) \geq 1-\frac{1}{n^p}\ \frac{1}{\left(\frac{\varepsilon}{C_p^{1/2p}}-\frac1{\sqrt{n}(\pee(A))^{1/2p}}\right)^{2p}}
\]
for all $\varepsilon > \frac{C_p^{1/2p}}{\sqrt{n}(\pee(A))^{1/2p}}$. 
\et

\bigskip

\bpr
We apply Theorem \ref{HigherMomentsMarkovchains} to $f=\bar{d}(\cdot,A)$, which is a function defined on $E^n$. 
It is easy to check that $\delta_i(f)\leq 1/n$, $i=1,\ldots,n$.
We first estimate $\E(f)$ by using \eqref{kloot}, which gives (using the fact that $f_{|A}=0$)
\[
\E(f)\leq \frac{C_p^{1/2p}}{\sqrt{n} (\pee(A))^{1/2p}}.
\]
Now we apply \eqref{polyci} with $t=\varepsilon$,
$$
\pee(f> \varepsilon) \leq \frac{C_p}{n^p}\ \frac{1}{\left(\varepsilon -\frac{C_p^{1/2p}}{\sqrt{n} (\pee(A))^{1/2p}}\right)^{2p}}.
$$
The result then easily follows.
\epr

\bigskip

\br
For $p=1$, the theorem holds under the assumption of Theorem \ref{DevroyeMarkovchains}; see Remark \ref{strounf}.
\er

As we saw in Section \ref{hoc}, we cannot have Gaussian bounds for certain house of cards processes, but only moment estimates up to a critical order.
In particular, this means that we cannot have a Gaussian measure concentration of Hamming neighborhoods.
But in that case we can apply the previous theorem and get polynomial measure concentration.

\bigskip

\noindent {\bf Acknowledgment}. The authors thank E. Verbitskiy for useful discussions on the house of cards process, and an anonymous referee for useful remarks. 



\begin{thebibliography}{99}


\bibitem{adam}
R. Adamczak, A tail inequality for suprema of unbounded empirical
processes with applications to Markov chains.
Electron. J. Prob. {\bf 13}, 1000-1034, (2008).

\bibitem{bg}
S. Bobkov and F. G\"otze.
Exponential integrability and transportation cost related to logarithmic Sobolev inequalities. J. Funct. Anal. {\bf 163} (1999), 1--28. 

\bibitem{bouch}
S. Boucheron, O. Bousquet, G. Lugosi and P. Massart,
Moment inequalities for functions of independent random
variables. Ann. Prob. {\bf 33}, 514-560, (2005).

\bibitem{bfg}
X. Bressaud, R. Fern\'andez, and A. Galves.
Decay of correlations for non-H\"olderian dynamics. A coupling approach. 
Electron. J. Probab. {\bf 4} (1999), no. 3 (19 pp.). 

\bibitem{burk1}
D.L. Burkholder.
Sharp inequalities for martingales and stochastic integrals.
Colloque Paul L\'evy sur les Processus Stochastiques (Palaiseau, 1987).
Ast\'erisque No. 157-158 (1988), 75--94. 

\bibitem{chatterjee}
S. Chatterjee.
Stein's method for concentration inequalities. 
Probab. Theory Related Fields {\bf 138} (2007), no. 1-2, 305--321. 

\bibitem{cckr}
J.-R. Chazottes, P. Collet, C. K\"ulske and F. Redig.
Concentration inequalities for random fields via coupling.
Probab. Theory Related Fields {\bf 137} (2007), no. 1-2, 201--225.

\bibitem{collet}
P. Collet.
Variance and exponential estimates via coupling. 
Bull. Braz. Math. Soc. {\bf 37} (2006), no. 4, 461--475. 


\bibitem{arnaud}
H. Djellout, A. Guillin and L. Wu.
Transportation cost-information inequalities and applications to random dynamical systems and diffusions. 
Ann. Probab. {\bf 32} (2004), no. 3B, 2702--2732. 

\bibitem{douc}
R. Douc, A. Guillin and E. Moulines.
Bounds on Regeneration Times and Limit Theorems for Subgeometric Markov Chains.
Annales Inst. H. Poincar\'e, to appear.

\bibitem{douc2}
R. Douc, G. Fort, E. Moulines and P. Soulier, Practical drift conditions
for subgeometric rates of convergence. Ann. Appl. Prob. {\bf 14}, 1353-1377, (2004).

\bibitem{douc3}
R. Douc, E. Moulines and P. Soulier,
Computable convergence rates for sub-geometric ergodic Markov chains.
Bernoulli {\bf 13}, 831-848  (2007).

\bibitem{fey}
A.\ Fey-den Boer, R.\ Meester, Ronald, C.\ Quant, and F. Redig.
A probabilistic approach to Zhang's sandpile model.  Comm. Math. Phys. {\bf 280}, 351--388,  (2008). 

\bibitem{goldstein}
S. Goldstein.
A note on specifications.
Z. Wahrsch. Verw. Gebiete,  {\bf 46}, 45-51 (1978/79).

\bibitem{k1}
L. Kontorovich and K. Ramanan.
Concentration Inequalities for Dependent Random Variables via the Martingale Method,
prepint (2007), to appear in Ann. Probab.
	
\bibitem{k2}
L. Kontorovich.
Obtaining Measure Concentration from Markov Contraction, preprint, 2007 (arXiv:0711.0987).

\bibitem{ledoux}
M. Ledoux.
{\em The concentration of measure phenomenon}, Mathematical Surveys and Monographs {\bf 89}.
American Mathematical Society, Providence R.I., 2001.


\bibitem{ligg}
T.M.\ Liggett, {\em Interacting particle systems.} Reprint of the 1985 original. Classics in Mathematics. Springer-Verlag, Berlin, 2005.


\bibitem{mac}
C. McDiarmid. 
On the method of bounded differences, in {\em Surveys in Combinatorics} 1989, Cambridge University Press, 
Cambridge (1989) 148--188.

\bibitem{marton-1}
K. Marton.
Bounding $\overline d$-distance by informational divergence: a method to prove measure concentration. 
Ann. Probab. {\bf 24} (1996), no. 2, 857--866. 

\bibitem{marton0}
K. Marton.
A measure concentration inequality for contracting Markov chains. 
Geom. Funct. Anal. {\bf 6} (1996), no. 3, 556--571. [Erratum:  Geom. Funct. Anal.  7  (1997),  no. 3, 609--613.]

\bibitem{marton1}
K. Marton.
Measure concentration for a class of random processes.
Probab. Theory Related Fields {\bf 110} (1998), no. 3, 427--439.

\bibitem{mao} Y.H. Mao, Convergence rates in strong ergodicity for Markov processes.  Stochastic Process. Appl. {\bf 116}, no. 12, 1964--1976,  (2006). 




\bibitem{rio}
E. Rio.
In\'egalit\'es de Hoeffding pour les fonctions lipschitziennes de suites d\'ependantes. 
[Hoeffding inequalities for Lipschitz functions of dependent sequences]
C. R. Acad. Sci. Paris S\'er. I Math. {\bf 330} (2000), no. 10, 905--908. 

\bibitem{samson}
P.-M. Samson.
Concentration of measure inequalities for Markov chains and $\Phi$-mixing processes.
Ann. Probab. {\bf 28} (2000), no. 1, 416--461. 

\bibitem{iceman}
H. Thorisson.
{\em Coupling, stationarity, and regeneration}.
Probability and its Applications (New York). Springer-Verlag, New York, 2000.

\end{thebibliography}
\end{document}